\documentclass[12pt]{amsart}
\usepackage{amssymb,verbatim}
\usepackage{epsfig}

\newtheorem{theorem}{Theorem}[section]
\newtheorem{corollary}[theorem]{Corollary}
\newtheorem{lemma}[theorem]{Lemma}
\newtheorem{claim}[theorem]{Claim}
\newtheorem{proposition}[theorem]{Proposition}

\newtheorem{example}[theorem]{Example}
\newtheorem{remark}[theorem]{Remark}

\theoremstyle{definition}
\newtheorem{definition}[theorem]{Definition}
\newtheorem{fact}[theorem]{Fact}

\newcommand{\To}{\rightarrow}

\newcommand{\orb}{\text{orb}}

\newcommand{\bbr}{\mbox{$\mathbb{R}$}}
\newcommand{\bbz}{\mbox{$\mathbb{Z}$}}

\newcommand{\bbn}{\mbox{$\mathbb{N}$}}

\newcommand{\di}{\text{d}}

\newcommand{\dP}{\delta P}
\newcommand{\dQ}{\delta Q}

\begin{document}

\title{Does a billiard orbit determine its (polygonal) table?}

\author{Jozef Bobok}

\author{Serge Troubetzkoy}

\address{KM FSv \v CVUT, Th\'akurova 7, 166 29 Praha 6, Czech Republic}
\email{bobok@mat.fsv.cvut.cz}

\address{Centre de physique th\'eorique\\
Federation de Recherches des Unites de Mathematique de Marseille\\
Institut de math\'ematiques de Luminy and\\
Universit\'e de la M\'editerran\'ee\\
Luminy, Case 907, F-13288 Marseille Cedex 9, France}
\email{troubetz@iml.univ-mrs.fr}
\urladdr{http://iml.univ-mrs.fr/{\lower.7ex\hbox{\~{}}}troubetz/} \date{}

\begin{abstract}
We introduce a new equivalence relation on the set of all polygonal
billiards. We say that two billiards (or polygons) are order
equivalent if each of the billiards has an orbit whose footpoints
are dense in the boundary and the two sequences of footpoints of
these orbits have the same combinatorial order. We study this
equivalence relation with additional regularity conditions on the
orbit.
\end{abstract}

\maketitle

\section{Introduction}

In mathematics one often wants to know if one can reconstruct a
object (often a geometric object) from certain discrete
data.  A famous example of this is the celebrated problem posed by Mark
Kac, ``Can one
hear the shape of a drum'', i.e., whether one can reconstruct a
drum head from knowing the frequencies at which it vibrates \cite{K}.
This problem was resolved negatively by Milnor in  dimension 16
\cite{M} and
then by Gordon, Webb and Wolpert
in dimension 2 \cite{GWW}.
Another well known example is a question posed by Burns and Katok
whether  a negatively
curved surface is determined by its marked length spectrum \cite{BK}.
The ``marked length spectrum'' of a surface $S$ is the function that
associates to each conjugacy class in $\pi_1(S)$ the length of the
geodesic in the associated free homotopy class. This question
was resolved positively by Otal \cite{O}.

In this article we ask if a polygonal billiard table is determined
by the combinatorial data of the footpoints of a billiard orbit. For
this purpose we introduce a new equivalence relation on the set of
all polygonal billiards.   Namely we say that two polygonal
billiards (polygons) are order equivalent if each of the billiards
has an orbit whose footpoints are dense in the boundary and the two
sequences of footpoints of these orbits have the same combinatorial
order. We study this equivalence relation with additional regularity
conditions on the orbit.
Our main results are the following, under a weak regularity
condition on the orbits, an irrational polygon can not be order
equivalent to a rational polygon and  two order equivalent rational
polygons must have the same number of sides with corresponding
corners having the same angle. In the case of triangles, only similar
triangles can be order equivalent.
In general, in the rational case,
one can not say more since any two rectangles are order equivalent.
If we furthermore assume that
the greatest common denominator of the rational angles is at least
$3$, then under a slightly stronger regularity condition we show that
the two order equivalent rational polygons must be similar.
In the case that the greatest common denominator is $2$, the two order equivalent rational
polygons must be affinely similar.

In Section 2 we start by recalling basic needed facts about
polygonal billiards while in Section 3 we prove various preparatory
lemmas. We prove our main results in Sections 4 (rational versus
irrational)
and 5 (rational versus rational).
Finally in Section 6 we summerize our results and ask some open questions.

\section{Polygonal Billiard}

A polygonal billiard table is a planar simply connected compact polygon $P$. The billiard
flow $\{T_t\}_{t\in\bbr}$ in $P$ is generated by the free motion of
a mass-point subject to the elastic reflection in the boundary. This
means that the point moves along a straight line in $P$ with a
constant speed until it hits the boundary. At a smooth boundary
point the billiard ball reflects according to the well known law of
geometrical optics: the angle of incidence equals to the angle of
reflection. If the billiard ball hits a corner, (a non-smooth
boundary point), its further motion is not defined. Additionally to
a corner, the billiard trajectory is not defined for a direction
tangent to a side.

By $D$ we denote
the group generated by the reflections in the lines through the
origin, parallel to the sides of the polygon $P$. The group $D$ is
either

\begin{itemize}\item finite, when  all the angles of $P$ are of the form $\pi
m_i/n_i$ with distinct co-prime integers $m_i$, $n_i$,  in this
case $D=D_N$ the dihedral group generated by the reflections in
lines through the origin that meet at angles $\pi/N$, where $N$
is the least common multiple of $n_i$'s,\end{itemize} or
\begin{itemize}\item countably infinite, when at least one angle
between sides of $P$ is an irrational multiple of $\pi$.
\end{itemize} In the two cases we will refer to the polygon as
rational respectively irrational.

Consider the phase space $P\times S^1$ of the billiard flow $T_t$,
and for $\theta\in S^1$, let $R_{\theta}$ be its subset of points
whose second coordinate belongs to the orbit of $\theta$ under $D$.
Since a trajectory changes its direction by an element of $D$ under
each reflection, $R_{\theta}$ is an invariant set of the billiard
flow $T_t$ in $P$.

The billiard map $T\colon~V_P=\cup e\times\theta\subset \dP\times
S^1\To V_P$ associated with the flow $T_t$ is the first return map
to the boundary $\dP$ of $P$. Here the union $\cup e \times \theta$
is taken over all sides of $P$ and for each side $e$ over the inner
pointing directions $\theta$
from $S^1$.  We will denote points of $V_P$ by $u = (x,\theta)$.

The set $P\times \theta$, resp.\ $\dP\times \theta$ will be called a
floor of the phase space of the flow $T_t$, resp.\ of the billiard
map $T$.

Let $\dP$ be oriented counterclockwise. For $x,x'\in \dP$, by
$[x,x']$ ($(x,x')$) we denote a closed (open) arc with outgoing
endpoint $x$ and incoming endpoint $x'$.

If $P,Q$ are simply connected polygons with counterclockwise
oriented boundaries, two sequences $\{x_n\}_{n\ge 0}\subset \dP$ and
$\{y_n\}_{n\ge 0}\subset \dQ$ have the same combinatorial order if
for each non-negative integers $k,l,m$
$$x_k\in [x_l,x_m]~\iff~y_k\in [y_l,y_m].$$

 The next definition introduces a new relation on the set of all simply connected
 polygons. The reader can verify  that it is reflexive, symmetric and
transitive, i.e., it is an equivalence relation. As usual, $\pi_1$
denotes the natural projection to the first coordinate (the foot
point).

\begin{definition}\label{d-relative}
We say that two polygons (or polygonal billiards) $P,Q$ are order
equivalent if there are points $u_0\in V_P$, $v_0\in V_Q$ such that
\begin{itemize}
   \item[(i)] $\overline{\{\pi_1(T^n(u_0))\}}_{n\ge 0}=\dP$, $\overline{\{\pi_1(S^n(v_0))\}}_{n\ge
   0}=\dQ$,
     \item[(ii)] the sequences $\{\pi_1(T^n(u_0))\}_{n\ge 0}$, $\{\pi_1(S^n(v_0))\}_{n\ge 0}$
  have the same combinatorial order.
  \end{itemize}
  In such a case we will write $(P;u_0)\approx (Q;v_0)$, the points $u_0,v_0$ will be sometimes called the leaders.
 \end{definition}

Let $t=\{x_n\}_{n\ge 0}$ be a sequence which is dense in $\dP$. The
$t$-address $a_t(x)$ of a point $x\in\dP$ is the set of  all increasing
sequences $\{n(k)\}_k$ of non-negative integers satisfying
$\lim_kx_{n(k)}=x$. It is clear that any $x\in\dP$ has a nonempty
$t$-address and $t$-addresses of two distinct points from $\dP$ are
disjoint.

For order equivalent polygons $P$, $Q$  with leaders $u_0,v_0$, we
will consider addresses with respect to the sequences given by
Definition \ref{d-relative}(ii):
$$t=\{\pi_1(T^n(u_0))\}_{n\ge 0},~s=\{\pi_1(S^n(v_0))\}_{n\ge 0}.$$

 It is an easy exercise to prove that the
map $\phi\colon\dP\to \dQ$ defined by
\begin{equation}\label{e3-projectphi}\phi(x)=y~\text{ if }a_t(x)=a_s(y)\end{equation}
is a homeomorphism.

Concerning the order-equivalence of polygonal billiards, the reader
can ask possible properties of corresponding equivalence classes. In
particular, can rational and irrational polygons simultaneously lie
in the same class? Can order equivalent polygons have different
(number of) angles, or at least different lengths of sides?

It follows immediately from the well known results that order
equivalent polygons need not be similar (see for example \cite{MT}).
\begin{example}\label{e-squarerectangle}Any two rectangles are order equivalent.\end{example}

We finish this section recalling several well known and useful (for
our purpose) results about polygonal billiards (see for example
\cite{MT}).
Recall that a flat
strip $\mathcal T$ is an invariant subset of the phase space of the
billiard flow/map such that
\begin{itemize}\item[1)] $\mathcal T$ is
contained in a finite number of floors,\item[2)] the billiard
flow/map dynamics on $\mathcal T$ is minimal in the sense that
any orbit which does not hit a corder is dense in $\mathcal T$,
\item[3)] the boundary of $\mathcal
T$ is non-empty and consists of a finite union of generalized
diagonals (a generalized diagonal is a billiard trajectory that
goes from a corner to a corner).\end{itemize}

By $\pi_2$ we denote the second natural projection (to the
direction). A direction, resp.\  a point $u$ from the phase space is
exceptional if it is the direction of a generalized diagonal, resp.\
$\pi_2(u)$ is such a direction. Obviously there are countably many
generalized diagonals hence also exceptional directions. A
direction, resp.\  a point $u$ from the phase space, which is not
exceptional will be called non-exceptional.

The set of the corners of $P$ is denoted by $C_P$.

\begin{proposition}\cite{MT}~\label{p-summary}Let $P$ be rational and $u_0\in V_P$. Then exactly one of the following three possibilities has to be satisfied.
\begin{itemize}\item[(i)] $u_0$ is periodic.
\item[(ii)] $\overline{\orb}(u_0)$ is a flat strip.
\item[(iii)] $\omega(u_0)=R_{\pi_2(u_0)}$; then
\begin{equation*}
\#(\{\pi_2(T^n(u_0))\colon~n\ge 0\})=2N,\end{equation*} and for
every $x\in\dP\setminus C_P$,
\begin{equation*}
\#\{u\in\omega(u_0)\colon~\pi_1(u)=x\}=N,\end{equation*} where $N$
is the least common multiple of the denominators of angles of
$P$. Moreover, in this case
\begin{equation*}\pi_2(\{u\in\omega(u_0)\colon~\pi_1(u)=x\})=\pi_2(\{u\in\omega(u_0)\colon~\pi_1(u)=x'\})\end{equation*}
whenever $x'\notin C_P$ belongs to the same side as $x$. The
possibility (iii) is true whenever $u_0$ is non-exceptional.
\end{itemize}
\end{proposition}

The billiard map $T$ has a natural invariant measure on its phase
space given by the formula $\mu=\sin\theta~\di x~\di\theta$, where
$\theta\in (0,\pi)$ is measured with respect to the direction of an
oriented side $e$. In a rational polygon we say that a point $u$ is
generic if it is non-exceptional, has bi-infinite orbit
and the billiard  map  restricted to the invariant surface $R_{\pi_2(u)}$
has a single invariant measure (this measure is then automatically the
measure $\mu$).

Finally we will use the following remarkable result from \cite{V}.
\begin{theorem}\cite{MT}~\label{p-veech}If $P$ is a regular $n$-gon
($n  \ne 3$) then any non-exceptional $u$
with infinite orbit is generic. Any
exceptional $u$, the trajectory of which does not hit a corner, is
periodic.\end{theorem}

\section{Preparatory Lemmas}

To define some important objects in dynamical systems, for example
an $\omega$-limit set - the corresponding phase space has to be
equipped with a topology providing  suitable convergence. In the
context of polygonal billiards it often means the simultaneous
convergence of foot points and directions (shortly pointwise
convergence). Since it simplifies several our proofs, throughout the
whole paper we will use the following modified definition of
convergence.

\begin{definition}\label{d-conver}Let $\{T^n(u_0)\}_{n\ge
0}$ be an infinite trajectory. The sequence $\{T^{n(k)}(u_0)\}_{k\ge
0}$ converges if and only if \begin{itemize}\item[(i)]
$\pi_1(T^{n(k)}(u_0))\to_k x\in\dP$,
\item[(ii)]
$\pi_1(T^{n(k)+1}(u_0))\to_k x'\in\dP$,
\item[(iii)] the points $x$ and $x'$ do not belong to the same side of $P$.
\end{itemize}
\end{definition}

In particular, the above definition omits the case when a sequence
$\{T^{n(k)}(u_0)\}_{k\ge 0}$ satisfies (i),(ii) but $x=x'$. Not to
rewrite known results in slightly modified versions we will use the
following
\begin{fact}\label{c}  Despite the fact that
this notion is stronger than the usual convergence we will use
classical results, such  as Proposition \ref{p-summary}(ii),(iii).
All the results we use remain true under our notion of
convergence.\end{fact}

\begin{lemma}\label{l-sides}Let $(P;u_0)\approx (Q;v_0)$, assume that to any point $x'\in
C_P\cup\phi^{-1}(C_Q)$ there exists a sequence $\{n(k)\}_{k\ge 1}$
such that the items (i),(ii) of Definition \ref{d-conver} are
satisfied, $x\neq x'$ and $x'$ is a both-side limit point of
$\{\pi_1(T^{n(k)+1}(u_0))\}_{k\ge 0}$. Then
$\phi(C_P)=C_Q$.\end{lemma}
\begin{proof}
Choose $x'\in C_P\setminus\phi^{-1}(C_Q)$. Let
$\{T^{n(k)}(u_0)\}_{k\ge 0}$ be a sequence satisfying the
assumptions of the lemma. Consider its `counterpart' in $Q$,
$\{S^{n(k)}(v_0)\}_{k\ge 0}$. Let
$$y=\phi(x)=\lim_k\pi_1(S^{n(k)}(v_0)),~y'=\phi(x')=\lim_k\pi_1(S^{n(k)+1}(v_0)).$$
Then  $\lim_k\pi_1(T^{n(k)+2}(u_0))$ does not exist since $x'$ is a
corner and a both-sided limit.  However $\lim_k\pi_1(S^{n(k)+2}(v_0))$ does exist,
this is impossible for order equivalent polygons $(P;u_0)\approx
(Q;v_0)$.

The case when $x'\in \phi^{-1}(C_Q)\setminus C_P$ can be excluded
analogously. Thus,  $C_P=\phi^{-1}(C_Q)$, i.e., $\phi(C_P)=C_Q$.
\end{proof}

As usual, an $\omega$-limit set of a point $u$  is denoted by
$\omega(u)$. Let us define the map
$\Phi\colon~\{T^n(u_0)\}_{n\ge 0}\to \{S^n(v_0)\}_{n\ge 0}$ by
\begin{equation}\label{e13-globalphi}\Phi(T^n(u_0))=S^n(v_0).\end{equation}

\begin{lemma}\label{l-sameconver} Let $(P;u_0)\approx (Q;v_0)$ and
$\phi(C_P)=C_Q$. The following hold
  \begin{itemize}\item[(i)]$T^{n(k)}(u_0)\To_k u\in V_P$ if and only if $S^{n(k)}(v_0)\To_k v\in V_Q$.
\item[(ii)] The map $\Phi$ from (\ref{e13-globalphi}) can be extended homeomorphically to the map
$\Phi\colon~\omega(u_0)\to \omega(v_0)$ satisfying (for all
$n\in\bbz$ for which the image is defined)
\begin{equation*}\label{e-conjugacy}\Phi(T^n(u))=S^n(\Phi(u)),~u\in\omega(u_0)
;\end{equation*} in particular, for every $x\in\dP$ we have
\begin{equation}\label{e-dihnumber}\Phi(\{u\in\omega(u_0)\colon~\pi_1(u)=x\})=\{v\in\omega(v_0)\colon~\pi_1(v)=\phi(x)\}.\end{equation}
 \item[(iii)] For every $u,u'\in\omega(u_0)$ with $\pi_1(u)=\pi_1(u')\notin C_P$,
\begin{equation}\label{e-order}\pi_2(u)<\pi_2(u')~\iff~\pi_2(\Phi(u))<\pi_2(\Phi(u')),\end{equation}
\item[(iv)] $u_0$ is recurrent if and only if $v_0$ is recurrent.
 \end{itemize}
 \end{lemma}
\begin{proof}(i) Let $T^{n(k)}(u_0)\To_k u$ and $x,x'$ be points from Definition \ref{d-conver}. Then both the limits
$y=\phi(x)=\lim_k\pi_1(S^{n(k)}(v_0))$ and
$y'=\phi(x')=\lim_k\pi_1(S^{n(k)+1}(v_0))$ exist and, since
$\phi(C_P)=C_Q$, the points $y,y'$ cannot lie on the same side of
$Q$. Thus  $S^{n(k)}(v_0)\To_k v$ where $v$ is the unit vector with foot
point $y$ and pointing towards $y'$.
The opposite implication is symmetric.

(ii) Denote $u_n=T^n(u_0)$. For $u\in\omega(u_0)$ for which
$u=\lim_ku_{n(k)}$ put
$$\Phi(u)=\lim_k\Phi(u_{n(k)});$$ using (\ref{e13-globalphi}) and
part (i) of the lemma one can easily show that $\Phi$ is well defined and a
homeomorphism. Thus (\ref{e-dihnumber}) follows.

(iii) Since the homeomorphism $\phi$ preserves the corners, we have
$\pi_1(\Phi(u))=\pi_1(\Phi(u'))\notin C_Q$. From (ii) we obtain that
$\pi_2(u)\neq \pi_2(u')$
 if and only if $\pi_2(\Phi(u))\neq \pi_2(\Phi(u'))$. If
 (\ref{e-order}) does not hold then there are positive integers
 $n,n'$ such that $u_n$, resp.\  $u_{n'}$ is close to $u$, resp.\  $u'$ (with the foot points on the same side, in particular) and either
\begin{equation*}\label{e-order1}\pi_2(u_n)<\pi_2(u_{n'})~\&~\pi_2(\Phi(u_n))>\pi_2(\Phi(u_{n'})),\end{equation*}
or an analogous possibility with changed signs of inequality is
true. Without loss of generality  consider the first possibility. Then for the
corresponding arcs
$$[\pi_1(u_{n'+1}),\pi_1(u))\cap (\pi_1(u),\pi_1(u_{n+1})]=\emptyset$$
and
$$[\pi_1(\Phi(u_{n'+1})),\pi_1(\Phi(u)))\cap (\pi_1(\Phi(u)),\pi_1(\Phi(u_{n+1}))]=$$ $$=[\pi_1(\Phi(u_{n'+1})),\pi_1(\Phi(u_{n+1}))]\neq\emptyset,$$
what is impossible for order equivalent polygons $(P;u_0)\approx
(Q;v_0)$.

Part (iv) easily follows from (i).
\end{proof}


Let $g$ be a function defined on a neighborhood of $y$. The derived
numbers $D^+g(y)$, $D_+g(y)$ of $g$ at $y$ are given by $$D^+g(y) =
\limsup_{h\to 0_+}\frac{g(y + h)-g(y)}{h},~D_+g(y) = \liminf_{h\to
0_+}\frac{g(y + h)-g(y)}{h}$$ and the analogous limits from the left
are denoted by $D^-g(y)$, $D_-g(y)$.

Let $(z,y)$ be the coordinates of $\bbr^2$ and let $p_{a,b} \subset
\bbr^2$ be the line with equation $y = a + z \tan b$. For short we
denote $p_{y_0,g(y_0)}$ by $p_{g(y_0)}$.
\begin{lemma}\label{l-infinitesimal}Let $g\colon~(c,d)\to (-\frac{\pi}{2},\frac{\pi}{2})$ be a continuous function, fix $C\subset (c,d)$ countable.
Assume that for some $y_0$ one of the four possibilities
\begin{equation}\label{e10}D^+g(y_0)>0,~
D_+g(y_0)<0,~D^-g(y_0)>0,~D_-g(y_0)<0\end{equation} is fulfilled.
Then there exists a sequence $\{y_n\}_{n\ge 1}\subset (c,d)\setminus
C$ such that $\lim_ny_n=y_0$ and the set of crossing points
$\{p_{g(y_0)}\cap p_{g(y_n)}\colon~n\ge 1\}$ is bounded in the
$\bbr^2$.
\end{lemma}
\begin{proof}We will prove the conclusion when $\epsilon=D^+g(y_0)>0$. The three remaining
possibilities can be shown analogously.

By our assumption there exists a decreasing sequence $\{y_n\}_{n\ge
1}\subset (c,d)\setminus C$ such that $\lim_ny_n=y_0$ and for each
$n$
\begin{equation*}\label{e6}g(y_n)>g_n=g(y_0)+\frac{\epsilon}{2}(y_n-y_0)>g(y_0);\end{equation*}
It means that the crossing point $p_{g(y_0)}\cap p_{g(y_n)}$ is
closer to the point $(0,y_0)$ than the crossing point
$A_n=p_{g(y_0)}\cap p_{y_n,g_n}$. The first coordinate of
$A_n=(\tilde z_n,\tilde y_n)$ satisfies
\begin{equation*}\label{e7}y_0+\tilde z_n\tan g(y_0)=y_n+\tilde z_n\tan g_n\end{equation*}
hence
\begin{equation*}\label{e8}\tilde z_n=\frac{y_0-y_n}{\tan[g(y_0)+\frac{\epsilon}{2}(y_n-y_0)]-\tan g(y_0)}.\end{equation*}
Since
\begin{equation*}\label{e9}\lim_n\tilde z_n=\frac{-2}{\epsilon}\cos^2g(y_0),\end{equation*}
both the sets $\{A_n\colon~n\ge 1\}$, $\{p_{g(y_0)}\cap
p_{g(y_n)}\colon~n\ge 1\}$ are bounded in the $\bbr^2$.
\end{proof}

\begin{remark}\label{r-derivednum}Note that if $g'(y_0)\neq 0$ at
some point $y_0$ then at least one of the four possibilities in
(\ref{e10}) has to be fulfilled.\end{remark}

To apply this lemma to billiards, we remind the reader that notion of
an unfolded billiard trajectory. Namely,
instead of reflecting it in a side of $P$ one may reflect $P$ in
this side and unfold the trajectory to a straight line.

\begin{lemma}\label{l-finitefloors}Let $(P;u_0)\approx (Q;v_0)$ with $P$ rational. Then the set of directions $$\{\pi_2(S^n(v_0))\colon~n\ge 0\}$$
along the trajectory of $v_0$ is finite.
\end{lemma}

\begin{remark}In Lemma \ref{l-finitefloors} we do not assume $u_0$ is non-exceptional. \end{remark}

\begin{proof} We know that $\phi(\dP)=\dQ$, where $\phi$ is a homeomorphism defined in
(\ref{e3-projectphi}). Under our assumption we cannot exclude the
possibility $\phi(C_P)\neq C_Q$. So, it will be convenient to
consider extended sets of corners, $C_P^*$, $C_Q^*$, defined as
$$C_P^*=\{x\in\dP\colon~x\in C_P\text{ or }\phi(x)\in
C_Q\},~C_Q^*=\phi(C_P^*),$$ and the sides of $P$, resp.\  $Q$ will be
taken with respect to $C_P^*$, resp.\  $C_Q^*$.

By Definition \ref{d-relative}(i) the first projection of the
forward trajectory of $u_0$ is dense in $\dP$, so in particular,
$u_0$ is not periodic. In such a case the trajectory of $u_0$ is
minimal either in a flat strip $\mathcal{T}$ or an invariant surface
$R_{\pi_2(u_0)}$ - see Proposition
\ref{p-summary}. Since the
conclusion can be verified for both possibilities in a similar way,
we will only treat the flat strip $\mathcal{T}$ case.

For an interval $(a,b)\subset \dP$ and an $\alpha\in S^1$ denote
\begin{equation}\label{e11}I:=\{n\in\bbn_{0}\colon~\pi_1(T^n(u_0))\in (a,b)~\&~\pi_2(T^n(u_0))=\alpha\}.\end{equation}
By Proposition \ref{p-summary} there are finitely many directions
along the trajectory of $u_0$ and, by our assumption, the foot
points of it are dense in $\dP$; passing to a subinterval and
choosing a suitable $\alpha$ we can assume that $(a,b)$, resp.\
$\phi((a,b))=(c,d)$ is a subinterval of a ``side" $e$ of $P$, resp.\
a ``side'' $f$ of $Q$ and the sequence $\{\pi_1(T^n(u_0))\}_{n\in
I}$ is dense in $(a,b)$. Obviously,
\begin{equation*}\tau=(a,b)\times \{\alpha\}\subset
\omega(u_0),~\sigma=\Phi(\tau)\subset \omega(v_0);\end{equation*}
there is a countable subset $\tau_0$ of $\tau$ such that each point
from $\tau\setminus \tau_0$ has a bi-infinite
trajectory (either the forward or backward trajectory starting from
any point of $\tau_0$ finishes in a corner from $C_P^*$).

Clearly,  the sequence $\{\pi_1(S^n(v_0))\}_{n\in I}$ is dense in
$(c,d)=\phi((a,b))$. Define a continuous function $g\colon~(c,d)\to
S^1$ by
\begin{equation}\label{e-def(f)}
g(y)=\pi_2(\Phi(\phi^{-1}(y)\times \{\alpha\})).\end{equation}  In
what follows we will show that since the polygons $P,Q$ are order
equivalent, the continuous function $g$ has to be constant. It is
sufficient to show that $g'(y_0)=0$ whenever $y_0\in (c,d)\setminus
C$, where $C=\pi_1(\Phi(\tau_0))$ is countable.
 Choose the origin of $S^1$ to be
the perpendicular direction of the side of $Q$ containing $(c,d)$
and  fix $y_0\in (c,d)\setminus C$; then for a sufficiently small
neighborhood $U(y_0)$ of $y_0$, $g(U(y_0))\subset
(-\frac{\pi}{2},\frac{\pi}{2})$.

For $y' \in U(y_0) \setminus C$ consider the unfolded (bi-infinite)
billiard trajectory of $(y',g(y'))$ under the billiard flow
$\{S_t\}_{t\in\bbr}$ in $Q$. Via unfolding, this trajectory corresponds to the line
$p_{g(y')}$ with the equation $y=y'+z\tan g(y')$.
\begin{claim}\label{cl}There is no sequence $\{y_n\}_{n\ge 1}\subset
(c,d)\setminus C$ such that $\lim_ny_n=y_0$ and the set of crossing
points $\{p_{g(y_0)}\cap p_{g(y_n)}\colon~n\ge 1\}$ is
bounded.\end{claim}\begin{proof}Denote  by $S^m(y)$ the copy of
$S^m((y,g(y)))$ on $p_{g(y)}$. To the contrary, let $\{y_n\}_{n\ge
1}\subset (c,d)\setminus C$ satisfy $\lim_ny_n=y_0$ and the set of
crossing points $\{p_{g(y_0)}\cap p_{g(y_n)}\colon~n\ge 1\}$ is
bounded. Without loss of generality we can assume that the crossing points lie in the
positive part of $p_{g(y_0)}$ and that there exists a value $m$ for
which all crossing points from $\{p_{g(y_0)}\cap
p_{g(y_n)}\colon~n\ge 1\}$ lie on the link
$(S^{m-1}(y_0),S^m(y_0))\subset p_{g(y_0)}$. Choose a point $y_n$
such that each two points
$$S^i(y_0),~ S^i(y_n),~i=0,\dots,m$$ lie on the
same side of $Q$ and the crossing point $p_{g(y_0)}\cap p_{g(y_n)}$
lies on the link $(S^{m-1}(y_n),S^m(y_n))\subset p_{g(y_n)}$. Such a
$y_n$ does exist since the function $g$ is continuous and
$p_{g(y_0)}$ does not contain any copy of a corner from $C_Q^*$.
Thus,
\begin{equation}\label{e-links}(S^{m-1}(y_0),S^m(y_0))\cap (S^{m-1}(y_n),S^m(y_n))\neq\emptyset\end{equation}
At the same time, for the points $x_0=\phi^{-1}(y_0)$,
$x_n=\phi^{-1}(y_n)$, $\lim_nx_n=x_0$, the iterates
$T^i((x_0,\alpha)),~ T^i((x_n,\alpha))$ lie on the same side of $P$,
$i=0,\dots,m$, and the unfolded both-side billiard trajectories of
$(x_0,\alpha)$, $(x_n,\alpha)$ under the billiard flow
$\{T_t\}_{t\in\bbr}$ in $P$ are parallel lines. In particular for
the links,
$$(T^{m-1}(x_0),T^m(x_0))\cap (T^{m-1}(x_n),T^m(x_n))=\emptyset,$$
what contradicts (\ref{e-links}) for order equivalent polygons
$(P;u_0)\approx (Q;v_0)$ and the claim is proved.
\end{proof}
Applying Lemma \ref{l-infinitesimal} and Claim \ref{cl} we obtain
that the function $g$ defined in (\ref{e-def(f)}) satisfies
$g'(y_0)=0$ for every $y_0\in (c,d)\setminus C$ hence as a
continuous function has to be constant. It means that
\begin{equation}\label{e-gaps}S^n(v_0)\in\sigma\text{ and }\pi_2(S^n(v_0))=g(y_0)\text{ whenever }n\in
I\end{equation} and $I$ is given by (\ref{e11}). We assumed that the
trajectory of $u_0$ is minimal in a flat strip $\mathcal{T}$. It is
known that in such a case the gaps in $I$ are bounded. This fact
together with (\ref{e-gaps}) imply that the set of directions
$\{\pi_2(S^n(v_0))\colon~n\ge 0\}$ along the trajectory of $v_0$ is
finite. This proves the lemma.
\end{proof}

\section{Rational versus Irrational}

As in Section 2, by $D$ we denote the group generated by the
reflections in the lines through the origin, parallel to the sides
of the polygon $P$. The group $D$ is (countable) infinite if and
only if the polygon $P$ is irrational \cite{MT}. In our case when
$P$ is simply connected it is equivalent to the condition that some
angle of $P$ is not $\pi$-rational.

\begin{theorem}\label{t-irrat}Let $P$ be irrational, $u$ a point from the
phase space. \begin{itemize}\item[(i)] If the orbit of $u$ does not
hit a corner and $\theta=\pi_2(u)$ is non-exceptional then
$\{\pi_2(T^n(u))\colon~n\ge 0\}$ is infinite.
\item[(ii)] If $u$ is not periodic, but visits only a finite
number of floors then ($u$ is recurrent and) $\overline{\orb}(u)$ is
a flat strip.
\end{itemize}
\end{theorem}
\begin{proof}
(i)~Fix $P$ and $\theta$ as in the statement of the lemma. Enumerate
the infinite orbit $D\theta$ by $\{\theta_1,\theta_2,\dots\}$.
Consider $I := \cup e_j \times \theta_i$. This union is taken over
all sides and for each side over the inner pointing (at some
interior point of the side) $\theta_i$ from $S^1$. Equip each $e_j
\times \theta_i$ by a measure $\sin\theta_i\, dx$ where $dx$ is the
arc length along $e_j$ and $\theta_i\in (0,\pi)$ is measured with
respect to the oriented side $e_j$. We will refer to this measure as
length. The billiard map $T$ leaves $I$ invariant and acts as
an infinite interval exchange map.  The boundary point of each
interval of continuity corresponds to a (first) preimage of a corner of $P$.

Suppose that the result is not true, then there is a point $x\in
\partial P$ such that the orbit of $(x,\theta)$ is infinite but
takes on only a finite number of angles $\{\theta_{i_1}, \dots,
\theta_{i_n}\}$. Let $I(x):= \cup e_j \times \theta_{i_k} \subset I$
be the (finite) set of intervals in the directions visited by the
orbit of $(x,\theta)$.

First note that the set $I(x)$ cannot be $T$-invariant, this would
contradict the fact that the $D$ orbit of $\theta$ is infinite.

Let $J^\pm := J^\pm(x) := \{ (x',\theta') \in I(x): \ T^\pm(x',\theta') \in
I(x)\}$. Each of these sets is a finite union of intervals and the
total length of $J^+$ is equal to the total length of $J^-$. Thus
the total length of the intervals of $I(x)$ where the forward map is
not defined is equal to the total length of the intervals where the
backwards map is not defined. Thus we can (in an arbitrary manner)
complete the definition of the partially defined map to an interval
exchange transformation $G$ (IET). The IET $G$
agrees with the partially defined first return billiard
map whenever it was defined.


Since $G$ is an IET, the well known topological
decomposition holds, the interval of definition is decomposed into
periodic and minimal components, with the boundary of the components
consisting of saddle connections, i.e., orbits starting and ending
at a point of discontinuity of the IET (see for example \cite{MT}).
Saddle connections correspond to generalized diagonals of the
billiard.

Since the original
billiard direction is non-exceptional the $G$-orbit of $(x,\theta)$
cannot be periodic.  By the topological decomposition theorem
the if orbit of $(x,\theta)$  is not all of
$I(x)$ it must accumulate on a saddle connection of $G$.
This saddle connection must be in the closure of the set $J^+$,
thus it corresponds to a generalized diagonal of the billiard map.
This contradicts the fact that $\theta$ is non-exceptional.

(ii) The proof is similar as for (i).  Define the ghost map and its
extension in the same way as above.  The conclusion is just the
consequence of the topological decomposition theorem for IETs.
\end{proof}

The main result of this section follows.
\begin{theorem}\label{t-ratirrat}Let $(P;u_0)\approx (Q;v_0)$, $P$ rational, $u_0$ non-exceptional.
Then $Q$ is rational with $v_0$ non-exceptional.
\end{theorem}\begin{proof}
All the assumptions of Lemma \ref{l-sides} are fulfilled, thus
$\phi(C_P)=C_Q$. By Lemma \ref{l-finitefloors}, the $v_0$ visits a
finite number of floors.

First let us assume that $v_0$ is exceptional, i.e., $v_0$ is parallel
to a generalized diagonal $d$. Let $y$, resp.\  $y'$ be an outgoing,
resp.\  incoming corner of $d$ with $y'=\pi_1(S^m((y,\beta)))$ for
some $m\in\bbn$ and a direction $\beta$ with respect to a side $f$.
Since by Lemma \ref{l-sameconver}(iv) also $v_0$ is recurrent and it
visits a finite number of floors, there is a sequence
$\{n(k)\}_{k\ge 0}$ such that for each $k$,
$(y_k,\beta)=S^{n(k)}(v_0)\in\omega(v_0)$ with $y\neq y_k\in f$ and
$\lim_ky_k=y$. Put $e=\phi^{-1}(f)$, $x_k=\phi^{-1}(y_k)$,
$x=\phi^{-1}(y)$ and $x'=\phi^{-1}(y')$. Then $x,x'\in C_P$ and
since also $u_0$ visits a finite number of floors, the sequence
$\{\pi_2(\Phi^{-1}((y_k,\beta)))\}_{k\ge 0}$ has to be eventually
constant, i.e.,
$\Phi^{-1}((y_k,\beta))=(x_k,\alpha)=T^{n(k)}(u_0)\in\omega(u_0)$
for each sufficiently large $k$. The reader can verify that $u_0$ is
parallel to a generalized diagonal outgoing from $x$ and incoming to
$x'$. This is impossible for $u_0$ non-exceptional. It shows that
$v_0$ has to be non-exceptional.

If the polygon $Q$ were irrational, Theorem \ref{t-irrat}(i) would
imply, since the orbit of $v_0$ is infinite, that $v_0$ is
exceptional. This is impossible by Theorem \ref{t-irrat}.
\end{proof}

\section{Rational versus Rational}

In a rational polygon a billiard trajectory may have only finitely
many different directions. In Section two we introduced the
invariant subset $R_{\theta}$ of the phase space
consisting of the set of
points whose second projection belongs to the orbit of $\theta$
under the dihedral group $D_N$; $R_{\theta}$ has the structure of a surface.
For non-exceptional $\theta$'s the
faces of $R_{\theta}$ can be glued according to the action of $D_N$
to obtain a flat surface depending only on the polygon $P$ but not
on the choice of $\theta$ - we will denote it $R_P$.

Let us recall the construction of $R_P$. Consider $2N$ disjoint
parallel copies $P_1,\dots,P_{2N}$ of $P$ in the plane. Orient the
even ones clockwise and the odd ones counterclockwise. We will glue
their sides together pairwise, according to the action of the group
$D_N$. Let $0<\theta=\theta_1<\pi/N$ be some angle, and let
$\theta_i$ be its $i$-th image under the action of $D_N$. Consider
$P_i$ and reflect the direction $\theta_i$ in one of its sides. The
reflected direction is $\theta_j$ for some $j$. Glue the chosen side
of $P_i$ to the identical side of $P_j$. After these gluings are
done for all the sides of all the polygons one obtains an oriented
compact surface $R_P$.

Let $p_i$ be the $i$-th vertex of $P$ with the angle $\pi m_i/n_i$
and denote by $G_i$ the subgroup of $D_N$ generated by the
reflections in the sides of $P$, adjacent to $p_i$. Then $G_i$
consists of $2n_i$ elements. According to the construction of $R_P$
the number of copies of $P$ that are glued together at $p_i$ equals
to the cardinality of the orbit of the test angle $\theta$ under the
group $G_i$, that is, equals $2n_i$.

 Each
corner $p\in C_P$ of $P$ corresponds to an angle $A(p)\in
(0,2\pi)\setminus \{\pi\}$.

\begin{proposition}\label{p-sameangles}Let $(P;u_0)\approx (Q;v_0)$, $P$ rational with $u_0$
non-exceptional. Then $\phi(C_P)=C_Q$ and $A(p)=A(\phi(p))$ for each
$p\in C_P$.
\end{proposition}

A triangle is determined (up to similarity) by its angles, thus
Proposition \ref{p-sameangles} implies
\begin{corollary}\label{c-triangle}
Let $(P;u_0)\approx (Q;v_0)$, $P$ a rational triangle with $u_0$
non-exceptional. Then $Q$ is similar to $P$.
\end{corollary}

\begin{proof} Theorem \ref{t-ratirrat} implies that also $Q$ is
rational with  a non-exceptional leader $v_0$. From Lemma
\ref{l-sides} follows that $\phi(C_P)=C_Q$; let $k=\# C_P=\#C_Q$,
number the corners $p_i$ of $P$, resp.\  $q_i$ of $Q$ to satisfy
$\phi(p_i)=q_i$, $i=1,\dots,k$. Since $P$, $Q$ are rational and
simply connected, $A(p_i)=\pi m_i^P/n_i^P$ and $A(q_i)=\pi
m_i^Q/n_i^Q$, where $m_i^P$, $n_i^P$, resp.\  $m_i^Q$, $n_i^Q$ are
coprime integers. In what follows, we will show that $n_i^P=n_i^Q$
and $m_i^P=m_i^Q$.

First of all, it is known that, since $u_0$ is non-exceptional, the
least common multiple $N_P$ of $n_i^P$'s is equal to
\begin{equation}\label{e-lcm}\#\{u\in\omega(u_0)\colon~\pi_1(u)=x\}=\#\{u\in\omega(u_0)\colon~\pi_1(u)=x'\},\end{equation}
whenever $x,x'\in\dP\setminus C_P$; moreover, if in addition $x$ and
$x'$ are from the same side then
\begin{equation}\label{e-lcmconstant}\pi_2(\{u\in\omega(u_0)\colon~\pi_1(u)=x\})
=\pi_2(\{u\in\omega(u_0)\colon~\pi_1(u)=x'\}),\end{equation} The
analogous equalities are true for $N_Q$.

Using (\ref{e-lcm}) and Lemma
\ref{l-sameconver}(ii)(\ref{e-dihnumber}) we get $N_P=N_Q=N$. Thus,
both rational billiards correspond to the same dihedral group $D_N$.

Second, consider the local picture around the $i$th vertex $p_i$.
Denote the two sides which meet at $p_i$ by $e$ and $e'$. Suppose
there are $2n^P_i$ copies of $P$ which are glued at $p_i$. Enumerate
them in a cyclic counterclockwise fashion $1,2,\dots,2n^P_i$. Since
$u_0$ is non-exceptional it's orbit is minimal, so it visits each of
the copies of $P$ glued at $p_i$. In particular the orbit crosses
each of the gluings (copy $j$ glued to copy $j+1$).

Now consider the orbit of $v_0$. We need to show that there are the
same number of copies of $Q$ glued at $q_i=\phi(p_i)$. Fix a
$j\in\{1,\dots,2n^P_i\}$ viewed as a cyclic group. Since $u_0$ is
non-exceptional the orbit of $u_0$ must pass from copy $j$ to copy
$j+1$ of $P$ or vice versa from copy $j+1$ to copy $j$. Suppose that
we are at the instant that the orbit $u_0$ passes from copy $j$ to
copy $j+1$ of $P$. At this same instant the orbit of $v_0$ passes
through a side. We label the two copies of $Q$ by $j$ and $j+1$
respectively. This labeling is consistent for each crossing from $j$
to $j+1$.

Since this is true for each $j$, the combinatorial data of the orbit
$u_0$ glue the corresponding $2n^P_i$ copies of $Q$ together in the
same cyclic manner as the corresponding copies of $P$. Note that the
common point of the copies of $Q$ is a common point of $\phi(e)$ and
$\phi(e')$ thus it is necessarily the point $q_i=\phi(p_i)$. In
particular, since Lemma \ref{l-sameconver}(iii) applies, we have
$2n^P_i$ copies of $Q$ glued around $q_i$ to obtain an angle which
is a multiple of $2\pi$. The total angle at this corner is by
definition $2\pi m^Q_i$. Thus $2n^Q_i$ must divide $2n^P_i$. The
argument is symmetric, thus we obtain $2n^P_i$ divides $2n^Q_i$. We
conclude that $n^P_i=n^Q_i$.

 Third, let us show that $m_i^Q= m^P_i$. Realizing the gluing of
 $2n_i^P$ copies of $P$ together at $p_i$ we get a point $p\in R_P$ with total angle of
 $2\pi m_i^P$. If $m_i^P>1$, the point $p$ is a cone angle
$2\pi m_i^P$ singularity. In any case, for the direction $\theta$
and the corresponding constant flow on $R_P$, there are
$m_i^P$ incoming trajectories
that enter $p$ on the surface $R_P$, hence also $m_i^P$
points in $V_P$ that finish their trajectory after the first iterate
at the corner $p_i$. Repeating all arguments for $Q$ and
$\vartheta=\pi_2(v_0)$, one obtain $m_i^Q$ points in $V_Q$ that
finish their trajectory after the first iterate at the corner
$q_i=\phi(p_i)$. Since such a number has to be preserved by the
homeomorphism $\Phi$, the inequality $m_i^P\neq m_i^Q$ contradicts
our assumption $(P;u_0)\approx (Q;v_0)$. Thus, $m_i^Q= m^P_i$.
\end{proof}

We have recalled in Section 2 that the billiard map $T$ has a
natural invariant measure on its phase space given by the formula
$\mu=\sin\theta~\di x~\di\theta$ and the notion of a generic point
$u$. In the case, when $P$ is rational and the corresponding
billiard flow is dense in the surface $R_P$,
 the measure $\mu$ sits on the
skeleton $K_P$ (union of edges) of $R_P$. In particular, an edge $e$
of $K_P$ associated with $\theta$ has the $\mu$-length $\vert
e\vert\cdot\sin\theta$. As before the number $N$ is defined as the
least common multiple of $n_i$'s, where the the angles of a simply
connected rational polygon $P$ are $\pi m_i/n_i$,

We have seen in Example \ref{e-squarerectangle} that for $N=2$ order-equivalent non-similar polygons exist 
even if they are uniquely ergodic - (note that rectangular billiards
are uniquely ergodic in every non excpetional direction). For any
rational polygon with $N=2$ we can speak - up to rotation - about 
horizontal, resp. vertical sides.  Two such polygons, $P$ and $Q$ with
sides $e_i$ resp.\ $f_i$, are  
{\em affinely similar}  if they have the same number of corners/sides, 
corresponding angles equal and there are positive numbers $a,b\in\bbr$ such that
$\vert e_i\vert/\vert f_i\vert=a$, 
resp.\ $\vert e_i\vert/\vert f_i\vert=b$ for any pair of corresponding horizontal, resp. vertical sides.

\begin{theorem}\label{t-similar}Let $(P;u_0)\approx (Q;v_0)$, $P$ rational, $u_0$ generic.
\begin{enumerate}
\item{}  If $N=N_P\ge 3$ then $Q$ is similar to $P$.
\item{}  If $N=N_P = 2$ the $Q$ is affinely similar to $P$. 
\end{enumerate}
\end{theorem}
\begin{proof}
1) We denote the sides of $P$, resp.\  $Q$ by $e_i$, resp.\  $f_i$. The
generic leader $u_0$ is, in particular, non-exceptional. Then
Theorem \ref{t-ratirrat} claims that also $Q$ is rational with $v_0$
non-exceptional and by Proposition \ref{p-sameangles},
$\phi(C_P)=C_Q$ and $A(p)=A(\phi(p))$ for each $p\in C_P$.

The billiard map $S$ has invariant measures defined on its phase
space; one of them is given by the formula $\nu=\sin\vartheta \di
f~\di\vartheta$. Consider the conjugacy $\Phi\colon~\omega(u_0)\to
\omega(v_0)$ guaranteed by Lemma \ref{l-sameconver}(ii). Since the
map $\Phi$ conjugates $T$ and $S$ and $\Phi(u_0)=v_0$,
$\nu=\Phi^*\mu$ and $\nu$ is in fact the unique invariant measure of
$S$ with $v_0$ generic.

Consider the reduction of the billiard map $T$ in the direction
$\pi_2(u_0)$ to an interval exchange. In our case when $u_0$ is
non-exceptional and, by Proposition \ref{p-summary}, the set of
directions $\Theta(e)=\pi_2(\pi^{-1}(x)\cap \omega(u_0))$, $x\in e$,
does not change on a side $e$, it consists in the union of $N$
copies of each side $e$ rescaled by the factor $\sin\theta_i$,
where $\theta_i\in\Theta(e)$ is measured with respect to the
direction of an oriented side $e$.  Let the interval $I$ be the
union of all rescaled interval-copies over all sides of $P$. The
billiard map $T$ induces a piecewise isometry $T_I$ of the segment
$I$ that can be modified \cite{MT} to an interval exchange
transformation. Using the fact that $T$ (measure $\mu$), resp.\  $S$
(measure $\nu$) is uniquely ergodic, the maps $T,S$ are conjugated
and lengths of interval-copies are in fact $\mu$-lengths of edges on
the skeleton $K_P$, resp.\  $\nu$-lengths of edges on the skeleton
$K_Q$ and the resulting interval exchange transformations are (after
an appropriate permutations of interval-copies in one of them) the
same. In particular, for a side $e$ of $P$, any $\theta\in\Theta(e)$
and $(f,\vartheta)=\Phi((e,\theta))$,
\begin{equation}\label{e-cosinus}\vert e\vert\sin\theta=\vert f\vert\sin\vartheta.\end{equation}
Assume that the least common multiple $N$ of the denominators of
angles of $P$ is greater than or equal to $3$. Since the maps $T,S$
are conjugated, the polygons $P$, $Q$ correspond to the same
dihedral group $D_N$ generated by the reflections in lines through
the origin that meet at angles $\pi/N$ The orbit of
$\theta^+_0=\pi_2(u_0)$, resp.\  $\vartheta^+_0=\pi_2(u_0)$ under
$D_N$ consists of $2N$ angles
$$\theta_j^+=\theta_0^++2j\pi/N,~\theta_j^-=\theta_0^-+2j\pi/N,$$
resp.\
$$\vartheta_j^+=\vartheta_0^++2j\pi/N,~\vartheta_j^-=\vartheta_0^-+2j\pi/N.$$
Since $N\ge 3$, for each side $e$, resp.\  $f$ one can find the
angles
$$\theta,\theta+2\pi/N\in\Theta(e),\text{ resp.\  }\vartheta,\vartheta+2\pi/N\in\Theta(f)$$ such
that by Lemma \ref{l-sameconver}(iii) $\Phi(e,\theta)=(f,\vartheta)$
and $\Phi(e,\theta+2\pi/N)=(f,\vartheta+2\pi/N)$. Then as in
(\ref{e-cosinus}),
\begin{equation*}\label{e-cosina}\vert e\vert\sin\theta=\vert f\vert\sin\vartheta,~\vert e\vert\sin(\theta+2\pi/N)=
\vert f\vert\sin(\vartheta+2\pi/N),\end{equation*} hence after some
routine computation we get $\vert e\vert=\vert f\vert$.

2) As in part 1) the
generic leader $u_0$ is non-exceptional, so by Theorem \ref{t-ratirrat} also $Q$ is rational with $v_0$
non-exceptional. By Proposition \ref{p-sameangles},
$\phi(C_P)=C_Q$ and $A(p)=A(\phi(p))$ for each $p\in C_P$. Similarly as in the mentioned proof, for a side $e$ of $P$, any $\theta\in\Theta(e)$
and $(f,\vartheta)=\Phi((e,\theta))=(\phi(e),\vartheta)$,
\begin{equation}\label{e-cosinus}\vert e\vert\sin\theta=\vert f\vert\sin\vartheta,\end{equation}
where $\theta$, resp. $\vartheta$ can be taken the same for any pair of horizontal, resp. vertical sides. Thus, the number $a=\vert e\vert/\vert f\vert$, resp.  $b=\vert e\vert/\vert f\vert$ does not depend on a concrete choice of a pair of $\phi$-corresponding horizonal , resp. vertical sides. This finishes the proof of our theorem.
\end{proof}
\section{Conclusions and open questions}

Let us summarize our main results in terms of order-equivalence
classes. We call an order equivalence class a non-exceptional
(resp.\ generic) order-equivalence class if the point $u_0$ in $P$ is
additionally supposed to be non-exceptional (resp.\ generic with
respect to $\mu$).
Two order equivalent polygons $P,Q$ are said to be
quasisimilar if $\phi(C_P)=C_Q$ and $A(p)=A(\phi(p))$ for each $p\in
C_P$.

\begin{theorem}Let us consider order-equivalence classes defined in Definition \ref{d-relative}. Then:\begin{itemize}
\item An irrational non-exceptional order-equivalence class does not contain rational polygons - Th.\ref{t-irrat}.
\item A rational non-exceptional order-equivalence class contains quasisimilar polygons - Prop.\ref{p-sameangles}.
\item A rational non-exceptional triangle
order-equivalence class is (up to the similarity) one point -
Cor.\ref{c-triangle}.
\item A rational generic
order-equivalence class with $N \ge 3$ is (up to the similarity) one
point - Th.\ref{t-similar}(1).
\item A rational generic
order-equivalence class with $N =2$ is (up to affine similarity) one
point - Th.\ref{t-similar}(2).
\item For each $n\ge 3$, ($n \ne 4$) the order-equivalence class of
  the regular $n$-gon is (up to the similarity) one
point.
\end{itemize}\end{theorem}

The last part of the theorem follows by combining Theorem
\ref{t-similar}
with  Proposition \ref{p-veech}.
The following questions remain open:

\begin{enumerate}
\item Can one replace the assumption that both orbits are dense
  (Definition \ref{d-relative} (i)) by only one orbit being
  dense?

\item Can one replace the assumption the orbit is generic by the orbit
being non-exceptional in Theorem \ref{t-similar}?

\item Can two irrational (non-similar) polygons be order equivalent?

\item For a rational polygon $P$, does the density of the (foot) orbit
in $\partial P$ imply that the  direction is non-exceptional?

\end{enumerate}

{\bf Acknowledgements}  We gratefully acknowledge the support of the
``poste tcheque'' of the University of Toulon. The first author was
also supported by MYES of the Czech Republic via contract MSM
6840770010.

\end{document}